\documentclass{article}
\usepackage{amsfonts, amsmath, amsthm}
\begin{document}
    \title{A Bijective Proof of a Major Index Theorem of Garsia and Gessel\footnote{This paper was written as a M.Sc. thesis under the direction of Gil Kalai of the Hebrew University and Yuval Roichman of Bar-Ilan University.  I extend my thanks to both of them for their time and assistance.}}
    \author{Moti Novick\\Department of Mathematics\\Hebrew University, Jerusalem}
    \date{\today}
    \maketitle
\section{Introduction}\label{sec:int}
In 1913, Percy MacMahon introduced the \emph{major index} statistic, defined for any permutation $\sigma=\sigma(1)...\sigma(n)$ of a multiset of integers of size $n$ as the sum of the descents of $\sigma$, i.e., $maj(\sigma)=\sum_{i=1}^{n-1}i\chi(\sigma(i)>\sigma(i+1))$ \footnote{We adopt the convention that for any statement $A$, $\chi(A)=1$ if $A$ is true and $\chi(A)=0$ if $A$ is false.}.  If $T$ denotes the multiset $\{1^{a_1},...k^{a_k}\}$ (i.e., the set of $a_i$ copies of the number $i$ for $i=1,...,k$, with $\sum_{i=1}^{k}a_i=n$), then MacMahon discovered that the generating function for the major index over the set $P(T)$ of permutations of $T$ is the following q-multinomial coefficient:
\begin{equation} \label{equ:mm} \sum_{\sigma \in P(T)} q^{maj(\sigma)} =
\left[
\begin{array}{c}
n \\ a_1,...,a_k
\end{array}
\right]
\end{equation}
He went on to prove in \cite{MacMahon} that this is also the generating function for the \emph{inversion number} of permutations of the same set (an \emph{inversion} of $\sigma$ is a pair $(i,j)\in[n-1]\times[n]$ with $i<j$ and $\sigma(i)>\sigma(j)$, and the inversion number is the total amount of such pairs: $inv(\sigma)= \sum_{i<j}\chi(\sigma(i)>\sigma(j))$).  This proved that these two statistics are equidistributed over all the permutations of any multiset of integers.  Specifically, in the case where $T=[n]:=\{1,...,n\}$ (i.e., $a_1=...=a_n=1$), and thus $P(T)=S_n$, we have:
\begin{equation} \label{equ:inv}
\sum_{\sigma\in S_n}q^{maj(\sigma)}= \sum_{\sigma\in S_n}q^{inv(\sigma)}= [n]_{q}! = (1+q)...(1+q+q^2+...+q^{n-1})
\end{equation}
Over fifty years passed before Foata \cite{Foata} discovered a bijective proof of this equidistribution result; a decade after that, he an Schutzenberger \cite{Foata_Schutz} showed that this bijection (when applied to $S_n$) preserves the inverse descent class of a permutation.  This proved that the two statistics are equidistributed over \emph{inverse descent classes} as well.
\medskip \newline
The following year, Garsia and Gessel stated the following theorem, an immediate consequence of
Stanley's theory of P-partitions \cite{Stanley2}.  MacMahon's results, as well as that of Foata and Schutzenberger, are corollaries of this theorem:
\medskip \newline
\textbf{Theorem} (\cite{Garsia_Gessel}, Theorem 3.1):  Let $\pi_1,...,\pi_k$ be ordered complementary subsets of $[n]$, where $\pi_i$ has length $a_i$ for $i=1,...,k$ (and hence $a_1+...+a_k=n$).  Let $\mathcal{S}(\pi_1,...,\pi_k)$ be the collection of permutations of $[n]$ obtained by shuffling $\pi_1,...,\pi_k$.  Then
\begin{equation}\label{equ:gg} \sum_{\sigma \in \mathcal{S}(\pi_1,...,\pi_k)}q^{maj(\sigma)} = \left[
\begin{array}{c}
n \\ a_1,...,a_k
\end{array}
\right]
q^{maj(\pi_1)+...+maj(\pi_k)}
\end{equation}
\medskip \newline
To see how MacMahon's results follow from this, define $a_0=0$, and let $\pi_i=(a_0+...+a_{i-1}, a_0+...+a_{i-1}+1,..., a_0+...+a_i)$.  Note that $maj(\pi_i)=0$ for all $i$ in this case, so the "$q$" term on the right is trivial.  MacMahon's results then follow by noting that there is a simple bijective correspondence between $\mathcal{S}(\pi_1,...,\pi_k)$ and the set of permutations of $\{1^{a_1},...k^{a_k}\}$ which preserves both inversion number and major index; just replace all elements of $\pi_i$ with the number $i$.
\medskip \newline
In this paper we provide a bijective proof of this theorem.  We actually prove the following theorem, from which the above result follows:
\medskip \newline
\textbf{Theorem 1.1}:  Let $\pi$ be an ordered subset of $[n]$ of length $a$, and let $\theta$ be an ordering of $[n]\backslash\pi$.  Let $\mathcal{S}(\theta,\pi)$ be the collection of permutations of $[n]$ obtained by shuffling $\theta$ and $\pi$.  Then
\begin{equation} \sum_{\sigma \in \mathcal{S}(\theta,\pi)}q^{maj(\sigma)} = \left[
\begin{array}{c}
n \\ a
\end{array}
\right]
q^{maj(\theta)+maj(\pi)}
\end{equation}
\medskip \newline
To obtain the version proved by Garsia and Gessel, simply apply Theorem 1.1 inductively on $i=2,...,k$ with $\theta$ some shuffle of $\pi_1,...,\pi_{i-1}$ and with $\pi=\pi_i$.  Summing over all $\theta$ and applying the inductive assumption then yields:
\begin{eqnarray}\sum_{\sigma \in \mathcal{S}(\pi_1,...,\pi_i)}q^{maj(\sigma)} &=& \left[
\begin{array}{c}
a_1+...+a_i \\ a_i
\end{array}
\right]
\left(
\left[
\begin{array}{c}
a_1+...+a_{i-1} \\ a_1,...,a_{i-1}
\end{array}
\right]
q^{maj(\pi_1)+...+maj(\pi_{i-1})}
\right)
q^{maj(\pi_i)}\nonumber\\
&=&
\left[
\begin{array}{c}
a_1+...+a_i \\ a_1,...,a_i
\end{array}
\right]
q^{maj(\pi_1)+...+maj(\pi_i)}\nonumber
\end{eqnarray}
\medskip \newline
When $i=k$ this clearly becomes equation (\ref{equ:gg}).
\medskip \newline
We will prove Theorem 1.1 by starting with the simplest case, namely that of $a=1$.  In this case the theorem essentially says that when any element $r\in[n]$ is inserted into a permutation $\sigma$ of $[n]\backslash r$, the resulting increase in major index is an element of $[n-1]_0:=\{0,1,...,n-1\}$ which depends uniquely on the index at which $r$ is inserted.  This is the simplest case to prove and leads immediately to a bijective proof of (\ref{equ:inv}).  We then proceed to the general case.  Finally, we return to the case of inverse descent classes, and show how our proof of the general case leads to a new bijective proof of that equidistribution result as well.
\section{Preliminaries}\label{sec:pre}
This section introduces the terminology and notation that will be used in the remainder of the paper (aside from what has been defined in the introduction).  An element $\sigma \in S_n$ is considered both as a word $a_1a_2...a_n$ (whose individual elements $a_1, a_2,...$ we call \emph{letters}), and as a bijection from $[n]$ to itself, with $\sigma(i)=a_i$ for $i = 1,...,n$.  A \emph{subword} of $\sigma$ is a string of distinct letters $a_1a_2...a_m$ for some $m \leq n$ such that $a_i \in [n]$ for $i=1,...,m$ and such that whenever $i<j\leq m$, $a_i$ precedes $a_j$ in $\sigma$.  The permutation $a_1a_2...a_n$ can be identified with the ordered sequence $\textbf{a}=(a_1,a_2,...,a_n)$, and conversely any ordered sequence of distinct integers can be identified with a permutation in the obvious manner.  The \emph{k-initial segment} of $\textbf{a}$ (for $k<n$) is the subsequence $(a_1,...,a_k)$.  $set(\textbf{a})$ denotes the (unordered) set of elements contained in $\textbf{a}$.
\medskip \newline
An index $i$ is a \emph{descent} of $\sigma$ if $a_i>a_{i+1}$, and the \emph{descent set} of $\sigma$ (denoted $Des(\sigma$)) is defined as $Des(\sigma) := \{i \in [n-1]: a_i>a_{i+1}\}$.  We denote by $d_k(\sigma)$ the number of descents in $\sigma$ greater than or equal to $k$ (i.e., the number of descents at or to the right of $\sigma(k)$; specifically, $d_1(\sigma)=|Des(\sigma)|$).  Indices of $\sigma$ that are not descents are called \emph{ascents}.  It is easily observed that $maj(\sigma)=\sum_{k=1}^{n}d_k(\sigma)$, as a descent at index $i$ is "accounted for" exactly $i$ times in the sum on the right.
\medskip \newline
We end this section by defining two new functions.  Our bijections will be based on a study of what happens to the major index of a permutation $\sigma$ of $[n]\backslash r$ when $r$ is inserted in the $k$-th position (i.e., before $\sigma(k)$, or at the right end if $k=n$) to create a new permutation which we denote $\sigma_k^r$.  Define $mi(\sigma, k, r) := maj(\sigma_k^r)-maj(\sigma)$ (the initials stand for \emph{major increment}).  We will also be interested in the \emph{major increment sequence} of $\sigma$ relative to $r$ defined as $MIS(\sigma,r) := (mi(\sigma,1,r),...,mi(\sigma,n,r))$.  In words, the major increment sequence of $\sigma$ relative to $r$ is the sequence of $n$ numbers whose $i$-th entry denotes the change in major index induced by inserting $r$ into $\sigma$ at the $i$-th position.
\begin{description}
\item[Example 2.1.] Inserting $r=7$ into the permutation $\sigma = 426351$, which has major index $9\,(=\,1+3+5)$.
\end{description}
\begin{center} \begin{tabular}{cccc}
$k$ & $\sigma_k^r$ & $maj(\sigma_k^r)$ & $mi(\sigma, k, r)$ \\
\hline
1 & 7426351 & 13 & 4 \\
2 & 4726351 & 12 & 3 \\
3 & 4276351 & 14 & 5 \\
4 & 4267351 & 11 & 2 \\
5 & 4263751 & 15 & 6 \\
6 & 4263571 & 10 & 1 \\
7 & 4263517 & 9  & 0
\end{tabular} \end{center}
Thus we have $MIS(\sigma,r)=(4,3,5,2,6,1,0)$.  Note that this sequence is a permutation of $[n-1]_0$.  This is no accident, as Proposition 3.1 will show.
\section{The Generating Function for Major Index Over $S_n$}
When $a=1$, we have \begin{equation}\left[ \begin{array}{c} n \\ a \end{array} \right]= \left[ \begin{array}{c} n \\ 1 \end{array} \right]= \frac{q^n-1}{q-1}=(1+q+...+q^{n-1}) \nonumber \end{equation} and thus Theorem 1.1 reads as follows: Given $n$ and $r\in[n]$, let $\theta$ be any permutation of $[n]\backslash r$.  Then
\begin{equation} \sum_{i=1}^{n}q^{maj(\theta_i^r)} = (1+q+...+q^{n-1})q^{maj(\theta)} \nonumber\end{equation}
Dividing both sides by $q^{maj(\theta)}$ yields:
\begin{equation} \sum_{i=1}^{n}q^{maj(\theta_i^r)-maj(\theta)} = \sum_{i=1}^{n}q^{mi(\theta, i,r)}=1+q+...+q^{n-1}\nonumber\end{equation}
In words:  The sequence $MIS(\theta,r) := (mi(\theta,1,r),...,mi(\theta,n,r))$ is a permutation of $[n-1]_0$.  This fact (phrased in very different terminology) was first noted by Gupta \cite{Gupta}.
\medskip \newline
Before proceeding with the proof, we explain how this special case gives a bijective proof of the equidistribution of inversion number and major index over $S_n$.  The \emph{inversion sequence} of $\sigma\in S_n$, denoted $I(\sigma)$, is the sequence whose $i$-th term $I(\sigma)(i)$ denotes the number of inversions of $\sigma$ whose first letter is $i$.  For example, if $\sigma=6257431$, then $I(\sigma)=(0,1,1,2,3,5,3)$.  It is immediately clear that $0\leq I(\sigma)(i)\leq i-1$ for $i=1,...,n$, and that $inv(\sigma)=\sum_{i=1}^{n}I(\sigma)(i)$.  Furthermore, any sequence $(a_1,...,a_n)$ with $0\leq a_i\leq i-1$ is the inversion sequence of some $\sigma\in S_n$-- construct $\sigma$ inductively by inserting the letter $i$ such that $a_i$ smaller letters lie to its right.  This immediately proves that $\sum_{\sigma\in S_n}inv(\sigma)= [n]_{q}! = (1)(1+q)...(1+q+q^2+...+q^{n-1})$, as every choice of exponents in the product on the right corresponds to an inversion sequence ($a_i$ coming from the $i$-th term in the product).
\medskip\newline
Suppose we are given $\sigma\in S_n$ with $I(\sigma)=(a_1,...,a_n)$.  Since $MIS(\pi,j)$ is a permutation of $[i-1]_0$ for all $i\leq n$, $\pi\in S_i$, and $0<j\leq i$, there is a unique permutation $\tau\in S_n$ which may be built by the successive insertion of the elements of $[n]$ \textbf{in any fixed order} such that the increase in major index at the $i$-th insertion is $a_i$.  Thus we have $inv(\sigma)=\sum_{i=1}^{n}a_i=maj(\tau)$.  As this mapping is easily reversible, it represents a bijection on $S_n$ mapping inversion number to major index, proving the equidistribution of these two statistics.
\begin{description}
\item[Example 3.1.] For $n=7$, we choose two orders of element insertion-- say, increasing order and the order $4-2-7-3-6-1-5$-- and we illustrate in both cases what permutation $\tau\in S_7$ corresponds to $\sigma=6257431$ (the example from above).  We construct $\tau$ so that, for $i=1,...,7$, the increase in major index resulting from the $i$-th insertion is the $i$-th element of $I(\sigma)=(0,1,1,2,3,5,3)$.
\end{description}
Insertion in increasing order yields:
\begin{center} \begin{tabular}{ll}
$1$ & $(maj=0)$\\
$21$ & $(maj=1)$\\
$231$ & $(maj=2)$\\
$4231$ & $(maj=4)$\\
$54231$ & $(maj=7)$\\
$542631$ & $(maj=12)$\\
$\tau=5472631$ & $(maj=15)$
\end{tabular} \end{center}
And insertion in the second order yields:
\begin{center} \begin{tabular}{ll}
$4$ & $(maj=0)$\\
$42$ & $(maj=1)$\\
$472$ & $(maj=2)$\\
$4372$ & $(maj=4)$\\
$64372$ & $(maj=7)$\\
$643721$ & $(maj=12)$\\
$\tau=6453721$ & $(maj=15)$
\end{tabular} \end{center}
Having shown what is implied by the following proposition, we now prove the proposition itself:
\medskip \newline
\textbf{Proposition 3.1}: If $r\in[n]$ and $\sigma$ is a permutation of $[n]\backslash r$, then $MIS(\sigma,r)$ is a permutation of $[n-1]_0$.
\medskip \newline
\textbf{Proof}:  We consider three cases: $r=n$, $r=1$, and $1<r<n$.  The first two cases will be proven directly, and the third will involve combining the first two.  Our proof is constructive in the sense that it not only proves that the sequence $MIS(\sigma,r)$ is a permutation of $[n-1]_0$ in each case but actually provides a method to construct this sequence.
\smallskip\newline
\underline{Case 1, $r=n$}:  Consider how $maj(\sigma)$ changes when $n$ is inserted at position $k$.  Clearly, for $k=n$ the change is zero, as $n > \sigma(n-1)$.  We consider the cases of $k<n$, i.e., the insertion of $n$ just before $\sigma(k)$. Firstly, the index of every element of $\sigma$ of index $k$ or higher increases by $1$ as a result of this insertion.  Specifically, every descent of $\sigma$ greater than or equal to $k$ increases by $1$ and hence the major index is increased by $d_k(\sigma)$.  In addition, $n$ itself (being larger than $\sigma(k)$ to its right) creates a new descent, but the consequent increment in $maj(\sigma)$ depends on whether $\sigma(k-1)<\sigma(k)$ or $\sigma(k-1)>\sigma(k)$.  If the former, or if $k=1$, then the insertion of $n$ produces a "brand new" descent at index $k$, increasing the major index by $k$.  If the latter, then the insertion produces a descent at index $k$ while eliminating a previously existing descent at index $k-1$ (as $\sigma_k^n(k-1)=\sigma(k-1)<n=\sigma_k^n(k)$).  Thus the additional increase in major index is $k-(k-1)=1$.  To summarize:
\begin{center}
\[
mi(\sigma,k,n)= \left\{
\begin{array}{llr}
0 & \mbox{if $k=n$.} & \\
d_k(\sigma)+k & \mbox{if $k-1$ is an ascent or $k=1$.} & (*)\\
d_k(\sigma)+1 & \mbox{if $k-1$ is a descent.} &
\end{array} \right.
\]
\end{center}
We show that $MIS(\sigma,n)$ is a permutation by constructing it explicitly backwards (i.e., from right to left).  The algorithm builds up $MIS(\sigma,n)$ using the decreasing index $i$ and intermediate permutations $\tau_n, \tau_{n-1},..., \tau_1$ (in that order) where each intermediate permutation is formed by appending a letter to the left of the previous one.  The algorithm, which we call \textbf{Algorithm L} for reasons to be explained later, proceeds as follows:
\begin{description}
\item[Step 1.] Let $A=1$, $B=n-1$, $i=n-1$, $\tau_n=0$.
\item[Step 2.] a) If $\sigma(i-1)>\sigma(i)$, define $\tau_i=A\tau_{i+1}$ and let $A = A+1$.\newline b) If $\sigma(i-1)<\sigma(i)$, or if $i=1$, define $\tau_i=B\tau_{i+1}$ and if $i\neq 1$ let $B = B-1$.
\item[Step 3.] Let $i=i-1$.
\item[Step 4.] If $i=0$, then output $\tau_1=MIS(\sigma,n)$; otherwise, return to Step 2.
\end{description}
We first note that the output of this algorithm is indeed a permutation of $[n-1]_0$.  The letters inserted over the $n-1$ iterations of the algorithm must be exactly the elements of $[n-1]$ because of the counter variables $A$ and $B$.  These start off as $1$ and $n-1$ respectively, then move toward each other at the rate of one "unit" per iteration (either $A$ increasing by $1$ or $B$ decreasing by $1$); thus after $n-2$ steps they are equal, and their common value becomes (in the iteration $i=1$) the first letter of $\tau_1$.  From this algorithm we see that $MIS(\sigma,n)$ is a permutation of $[n-1]_0$ with a special form: For any $i\in[n]$, there exist $A(i),B(i)\in[n-1]_0$ such that the first $i$ entries in $MIS(\sigma,n)$ are the numbers between $A(i)$ and $B(i)$ inclusive.  We call such a permutation an \textbf{A-B permutation}.
\medskip \newline
It remains to show that indeed $\tau_1=MIS(\sigma,n)$ as claimed.  The last element of $MIS(\sigma,n)$ is $mi(\sigma,n,n)=0=\tau_n$, by (*); thus these two sequences agree in their rightmost elements, i.e., after zero iterations of the algorithm.  We proceed by induction on the number of iterations.  Suppose that after $j$ iterations, $\tau_{n-j}$ is identical to the last $j+1$ elements of $MIS(\sigma,n)$.  The $(j+1)$-th letter from the end of $MIS(\sigma)$ is $mi(\sigma,n-j,n)$.  Suppose $n-j-1$ is a descent, so $mi(\sigma,n-j,n)=d_{n-j}(\sigma)+1$ (by (*)).  If $n-j-1$ is the rightmost descent in $\sigma$, then clearly $A=1$, and additionally $d_{n-j}(\sigma)=0$ so that $mi(\sigma,n-j,n)=1$, as desired.  If $n-j-1$ is not the rightmost descent in $\sigma$, suppose that the first descent to the right of index $n-j-1$ is at index $l$.  Then $A=(d_{l+1}(\sigma)+1)+1$, by the inductive assumption.  But the sum in parentheses is just $d_{n-j}(\sigma)$, since by assumption $n-j-1$ and $l$ are consecutive descents.  Thus $A=d_{n-j}(\sigma)+1=mi(\sigma,n-j,n)$, as we wished to prove.
\medskip \newline
The other possible case, namely that $n-j-1$ is an ascent or that $j=n-1$ (the latter meaning that the algorithm is up to the iteration $i=1$, the last one), proceeds in a parallel manner.  In this case, $mi(\sigma,n-j,n)=d_{n-j}(\sigma)+(n-j)$.  If $n-j-1$ is the rightmost ascent in $\sigma$ (or, if $j=n-1$ and $\sigma$ has no ascents), then $B=n-1$, and additionally $d_{n-j}(\sigma)=j-1$ so that $mi(\sigma,n-j)=(j-1)+(n-j)=n-1$, as desired.  If $n-j-1$ is not the rightmost ascent in $\sigma$ (or if $j=n-1$ and $\sigma$ has at least one ascent), suppose that the first ascent to the right of index $n-j-1$ is at index $l$.  Then $B=(d_{l+1}(\sigma)+(l+1))-1=d_{l+1}(\sigma)+l$, by induction.  We have $d_{n-j}(\sigma)=d_{l+1}(\sigma)+(l-n+j)$, because by assumption every one of the $l-n+j$ indices between $n-j-1$ and $l$ is a descent.  Thus $B=d_{l+1}(\sigma)+l=d_{n-j}(\sigma)+(n-j)$, as we wished to prove.  This completes the case of $j=n$.
\medskip \newline
Based on Algorithm L, we define for any permutation $\sigma$ of length $n$ and any $k\in[n]$, the ordered pair $L(\sigma,k)=(d_k(\sigma)+1,d_k(\sigma)+k)$.  $L(\sigma,k)$ is precisely the ordered pair of the counter variables $(A,B)$ at the start of the iteration that computes $mi(\sigma,k,n)$, which must be equal to one of these variables.  As noted, the first $k$ elements of $MIS(\sigma,n)$ are precisely the interval from one end of $L(\sigma,k)$ to the other.
\smallskip\newline
\underline{Case 2, $r=1$}:  The proof of this case proceeds in a similar manner. Consider how $maj(\sigma)$ changes when $1$ is inserted at position $k$.  Clearly, for $k=n$ the change is $n-1$, as $\sigma(n-1)>1$ and thus the index $n-1$ becomes a descent.  We consider the cases of $k<n$.  As before, the major index is increased by $d_k(\sigma)$ due to the increased index of every later descent.  In addition, inserting $1$ at index $k>1$ makes $k-1$ into a descent (as $1$ is smaller than all letters of $\sigma$), but the consequent increment in $maj(\sigma)$ depends on whether $\sigma(k-1)<\sigma(k)$ or $\sigma(k-1)>\sigma(k)$.  If the former, then the additional increase in $maj(\sigma)$ is $k-1$.  If the latter, then there is no such increase, as $k-1$ was already a descent in $\sigma$.  Also there is clearly no additional increase if $k=1$.  To summarize:
\begin{center}
\[
mi(\sigma,k,1)= \left\{
\begin{array}{llr}
n-1 & \mbox{if $k=n$.} & \\
d_k(\sigma)+(k-1) & \mbox{if $k-1$ is an ascent.} & (**) \\
d_k(\sigma) & \mbox{if $k-1$ is a descent or $k=1$.} &
\end{array} \right.
\]
\end{center}
We show that $MIS(\sigma, 1)$ is an A-B permutation of $[n-1]_0$ by demonstrating that it can be constructed from last element to first by an algorithm similar to Algorithm L.  Consider the following algorithm, named \textbf{Algorithm G} for reasons which will soon become clear:
\begin{description}
\item[Step 1.] Let $A=0$, $B=n-2$, $i=n-1$, $\tau_n=n-1$.
\item[Step 2.] a) If $\sigma(i-1)>\sigma(i)$ or if $i=1$, define $\tau_i=A\tau_{i+1}$ and if $i\neq 1$ let $A = A+1$. \newline b) If $\sigma(i-1)<\sigma(i)$, define $\tau_i=B\tau_{i+1}$ and let $B = B-1$.
\item[Step 3.] Let $i=i-1$.
\item[Step 4.] If $i=0$, then output $\tau_1=MIS(\sigma, 1)$; otherwise, return to Step 2.
\end{description}
The output $\tau_1$ of Algorithm G must be an A-B permutation of $[n-1]_0$ for the same reasons given regarding Algorithm L.  To see that $\tau_1= MIS(\sigma, 1)$, compare these two sequences starting from the right.  The last element of $MIS(\sigma,1)$ is $mi(\sigma,n,1)$, which by the first paragraph of this proof is $n-1= \tau_n$.  Thus the two sequences agree in their rightmost elements.
\medskip \newline
For the remaining elements, we could proceed  by induction as we did in the $r=n$ case, but there is a simpler proof.  The only difference between Algorithms L and G regards the initial conditions: In Algorithm L the initial values of the counter variables are $A=1$ and $B=n-1$, each greater by one than the corresponding initial value in Algorithm G (the case of $i=1$ is treated identically in both algorithms because $A=B$ by that iteration, as pointed out above).  A comparison of (*) and (**) reveals that $mi(\sigma,k,n)=mi(\sigma,k,1)+1$ for all $k<n$.  The fact that Algorithm L yields $MIS(\sigma,n)$ thus implies that Algorithm G yields $MIS(\sigma,1)$, as the case of $r=1$ is simply a "shift by $-1$" of the case of $r=n$, both in values of the counter variables and in corresponding major increments.  This completes the proof of the case $r=1$.
\medskip \newline
We define for any permutation $\sigma$ of length $n$ and any $k\in[n]$ the ordered pair $G(\sigma,k)=(d_k(\sigma),d_k(\sigma)+(k-1))$.  $G(\sigma,k)$ is precisely the ordered pair of the counter variables $(A,B)$ at the start of the iteration that computes $mi(\sigma,k,1)$.  We note that for any $\sigma$ and $k$, $L(\sigma,k)=G(\sigma,k)+(1,1)$.
\smallskip\newline
\underline{Case 3, $1<r<n$}: For the general case of $1<r<n$, we begin by partitioning $\sigma$ into maximal segments such that within each segment either every letter is less than $r$ (and the segment is denoted a \emph{lesser segment}) or every letter is greater than $r$ (and the segment is a \emph{greater segment}).  Since these segments are maximal, their order of appearance in $\sigma$ alternates between lesser and greater.
\begin{description}
\item[Example 3.2.] $\sigma=1762834$, $r=5$.  Then $\sigma$ is partitioned as $1\mid 76\mid 2\mid 8\mid 34$.
\item[Example 3.3.] $\sigma=5768312$, $r=4$.  Then $\sigma$ is partitioned as $5768\mid 312$.
\end{description}
The key observation is that within each lesser segment, the effect of inserting $r$ is the same as that of inserting $n$, and that within each greater segment the effect of inserting $r$ is the same as that of inserting $1$.  This explains the names of the algorithms: Algorithm L yields the major increment sequence when a letter is inserted into a lesser segment, while Algorithm G does the same when a letter is inserted into a greater segment.  Thus we expect that the appropriate algorithm in this case is one which "alternates" between Algorithms L and G in an appropriate manner; incredibly enough, the correct algorithm (which we dub the \textbf{L-G Algorithm}) does so while itself producing an $A-B$ permutation (we assume that from the outset $\sigma$ is partitioned into greater and lesser segments based on $r$):
\begin{description}
\item[Step 1.] Let $i=n-1$.
\item[Step 2.] \textbf{a)} If $\sigma(n-1)<r$, let $A=1$, $B=n-1$, $\tau_n=0$. \newline \textbf{b)} If $\sigma(n-1)>r$, let $A=0$, $B=n-2$, $\tau_n=n-1$.
\item[Step 3.] \textbf{a)} If $\sigma(i)<r$ and $\sigma(i-1)<r$, \textbf{or} if $i=1$ and $\sigma(i)<r$, perform step 2 of Algorithm L. \newline \textbf{b)} If $\sigma(i)>r$ and $\sigma(i-1)>r$, \textbf{or} if $i=1$ and $\sigma(i)>r$, perform step 2 of Algorithm G. \newline \textbf{c)} If $\sigma(i)>r$ and $\sigma(i-1)<r$, define $\tau_i=A\tau_{i+1}$ and let $A = A+1$. \newline \textbf{d)} If $\sigma(i)<r$ and $\sigma(i-1)>r$, define $\tau_i=B\tau_{i+1}$ and let $B = B-1$.
\item[Step 4.] Let $i=i-1$.
\item[Step 5.] If $i=0$, then output $\tau_1=MIS(\sigma, r)$; otherwise, return to Step 3.
\end{description}
The proof that this algorithm works is somewhat involved and technical, and we include it as an appendix.  Here we explain the method of the algorithm intuitively and mention an important property which will be needed later.  Step 2 sets up the initial conditions to be those of Algorithm L (if the rightmost segment of $\sigma$ is lesser) or those of Algorithm G (if that segment is greater).  Step 3 is the crucial one, evaluating the appropriate major increment depending on where $r$ is inserted.  It uses Algorithm L to insert $r$ within a lesser segment and Algorithm G to insert $r$ within a greater segment; this is the function of parts (a) and (b).  Parts (c) and (d) cover the "transition" steps of inserting $r$ between one kind of segment and the other.  As shown in the proof of the algorithm, the values of $A$ and $B$ after step 3 depend on what part of this step was implemented: If part (a) or (c) was implemented in the iteration $i=k+1$, then step 3 concludes with $(A,B)=L(\sigma,k)$; if part (b) or (d) was implemented, then step 3 concludes with $(A,B)=G(\sigma,k)$.  This fact will be vital in proving Lemma 4.1.
\section{The General Case}\label{sec:equinv}
We now use the results of the previous section to prove the general case (i.e., the case of $a>1$) of Theorem 1.1.  Let $b=n-a$ denote the length of $\theta$.  We will establish a bijection $\Phi$ between the set $\mathcal{S}(\theta,\pi)$ of shuffles of $\theta$ and $\pi$ and the set $\mathcal{P}(b,a)$ of partitions containing $a$ parts (some of which may be zero) all less than or equal to $b$.  Both of these sets have cardinality $\binom{a+b}{a}=\binom{n}{a}$ (in the case of $\mathcal{S}(\theta,\pi)$, an shuffle is determined uniquely by a choice of $a$ indices at which $\pi$ is inserted, and conversely; in the case of $\mathcal{P}(b,a)$, a partition $\lambda=(\lambda_1,...,\lambda_a)$ with $0\leq\lambda_1\leq...\leq\lambda_a\leq b$ becomes a uniquely determined $a$-element subset of $[n]$ by adding $i$ to $\lambda_i$, for $i=1,...,n$, and conversely).  Given $\lambda=(\lambda_1,...,\lambda_a)\in\mathcal{P}(b,a)$, denote the sum $\sum_{i=1}^{a}\lambda_i$ as $|\lambda|$.  Our bijection $\Phi: \mathcal{S}(\theta,\pi)\rightarrow \mathcal{P}(b,a)$ will have the property that for $\sigma\in\mathcal{S}(\theta,\pi)$,
\begin{equation}\label{equ:bijection}
maj(\sigma)=maj(\theta)+maj(\pi)+|\Phi(\sigma)|
\end{equation}
Raising $q$ to both sides of this equation, summing the left side over $\mathcal{S}(\theta,\pi)$ and the right side over $\mathcal{P}(b,a)$ (which preserves equality, by the bijection), yields:
\begin{equation}
\sum_{\sigma \in \mathcal{S}(\theta,\pi)}q^{maj(\sigma)} =
q^{maj(\theta)+maj(\pi)}\sum_{\lambda\in\mathcal{P}(b,a)}q^{|\lambda|}
\nonumber
\end{equation}
As is well-known, the generating function for the sums of the partitions in $\mathcal{P}(b,a)$ can be expressed as a $q$-binomial coefficient:
\begin{equation} \sum_{\lambda\in\mathcal{P}(b,a)}q^{|\lambda|} = \left[
\begin{array}{c}
b+a \\ a
\end{array}
\right]=
\left[
\begin{array}{c}
n \\ a
\end{array}
\right]\nonumber
\end{equation}
In fact, some sources define the $q$-binomial coefficient in this manner; see, e.g., (\cite{Bressoud}, chapter 3).  Thus, this bijection proves Theorem 1.1.
\medskip \newline
To define our bijection we need some new notation.  Given $\sigma\in\mathcal{S}(\theta,\pi)$, we imagine that $\sigma$ is constructed by the insertion of $\pi$ into $\theta$ one letter at a time, the letters being inserted in the \emph{reverse} of their order of appearance in $\pi$ (i.e., if $\pi=(\pi(1),...,\pi(a))$, then $\pi(a)$ is inserted first and $\pi(1)$ is inserted last).  Note that every insertion occurs to the left of the previous one.  Let $\sigma_i$ denote the subword of $\sigma$ consisting of $\theta$ and the elements $\pi(i),...,\pi(a)$, so that $\sigma_a,\sigma_{a-1},...,\sigma_1=\sigma$ represent the intermediate steps of the insertion procedure just described (as a convention, define $\sigma_{a+1}:=\theta$).  Let $k_i$ denote the position at which $\pi(i)$ is inserted into $\sigma_{i+1}$ to yield $\sigma_i$.  Since every insertion occurs to the left of the previous one, we have $k_1\leq...\leq k_a$.
\medskip\newline
With this construction procedure, let $m_i=maj(\sigma_i)-maj(\sigma_{i+1})$ (i.e., $m_i$ denotes the increase in major index induced by the insertion of $\pi(i)$) and let $t_i=m_i-d_i(\pi)$ .  We claim the following:
\medskip\newline
\textbf{Theorem 4.1}: The mapping $\Phi: \mathcal{S}(\theta,\pi)\rightarrow\mathcal{P}(b,a)$ defined by $\Phi(\sigma)=set((t_1,...,t_a))$ is a bijection between $\mathcal{S}(\theta,\pi)$ and $\mathcal{P}(b,a)$.
\begin{description}
\item[Example 4.1.]  Let $\theta=5274, \pi=631$, and $\sigma=5276341$.  Then $maj(\theta)=4$, and we have:
\end{description}
\begin{center} \begin{tabular}{llll}
$\sigma_3=52741$ & $k_3=5$ & $maj(\sigma_3)=8$ & $m_3=4$\\
$\sigma_2=527341$ & $k_2=4$ & $maj(\sigma_2)=9$ & $m_2=1$\\
$\sigma_1=5276341$ & $k_1=4$ & $maj(\sigma_1)=14$ & $m_1=5$\\
\end{tabular} \end{center}
Thus in this example, $t_1=5-2=3, t_2=1-1=0,$ and $t_3=4-0=4$, so $\Phi(\sigma)=\{0,3,4\}\in\mathcal{P}(b,a)$.
\medskip\newline
$\Phi$ indeed satisfies property (\ref{equ:bijection}):
\begin{equation}
maj(\sigma)-maj(\theta)= \sum_{i=1}^{a}m_i= \sum_{i=1}^{a}d_i(\pi)+\sum_{i=1}^{a}(m_i-d_i(\pi))=maj(\pi)+|\Phi(\sigma)| \nonumber
\end{equation}
It remains to show that $\Phi$ is indeed a bijection, and the remainder of this section is devoted to proving this fact.  It is not immediately clear that $\Phi$ even maps $\mathcal{S}(\theta,\pi)$ into $\mathcal{P}(b,a)$ at all.  To prove both that $\Phi(\sigma)\in\mathcal{P}(b,a)$ for all $\sigma\in\mathcal{S}(\theta,\pi)$ and that $\Phi$ is a bijection we need to have some idea of what the sequences $MIS(\sigma_{i+1},\pi_i)$ look like (for $i=1,...,a$).  To this end we prove the following lemma:
\medskip\newline
\textbf{Lemma 4.1}:  Let $\tau$ be a permutation of length $n-1$, $p,q\notin\tau$.  Let $\tau_j^p$ denote the permutation of length $n$ formed by the insertion of $p$ into $\tau$ at index $j$.  Then the first $j$ elements of $MIS(\tau_j^p,q)$ are some permutation of the set $\{x+\chi(q>p) | x$ is in the $j$-initial segment of $MIS(\tau,p)\}$.
\begin{description}
\item[Example 4.2.]  Let $\tau=436152$, $p=8$, $q=7$, and $j=5$.  Then $\tau_j^p=4361852$, and $\chi(q>p)=0$.  A quick calculation yields $MIS(\tau,8)=(4,3,5,2,6,1,0)$ and $MIS(\tau_j^p,7)=(5,4,6,3,2,7,1,0)$.  The first five elements of these two sequences are indeed the same.  Reversing the values of $p$ and $q$, we have $\tau_j^p=4361752$, $\chi(q>p)=1$, $MIS(\tau,7)=(4,3,5,2,6,1,0)$ and $MIS(\tau_j^p,8)=(5,4,6,3,7,2,1,0)$.  The first five elements of the first sequence, each increased by 1, yield the first five elements of the second sequence (coincidentally the order is preserved, but this need not be true in general).
\end{description}
\textbf{Proof (of lemma)}:  The first $j$ elements of $MIS(\tau,p)$ (resp., $MIS(\tau_j^p,q)$) comprise the subset of $[n-1]$ between the counter variables $A$ and $B$ (inclusive) after the iteration of the L-G algorithm which computes the $j+1$-th element of this sequence, namely $mi(\tau,j+1,p)$ (resp., $mi(\tau_j^p,j+1,q)$).  First consider the case of $p>q$, so that $\chi(q>p)=0$.  We wish to show that the respective values of $A$ and $B$ are the same at these points in the construction of the two sequences.  The letter $\tau_j^p(j)=p$ is part of a greater segment relative to $q$.  Thus the L-G algorithm uses either part (b) or (d) in Step 3 when computing $mi(\tau_j^p,j+1,q)$, and hence ends this iteration with $(A,B)=G(\tau_j^p,j)=(d_j(\tau_j^p),d_j(\tau_j^p)+j-1) $.  We wish to prove that this last equation holds after the L-G algorithm computes the element $mi(\tau,j+1,p)$ of $MIS(\tau,p)$ as well.  If $\tau(j)<p$, i.e., $\tau(j)$ is part of a lesser segment relative to $p$, then $mi(\tau,j+1,p)$ is computed using part (a) or (c) of Step 3, and so after this step $(A,B)=L(\tau,j)=(d_j(\tau)+1,d_j(\tau)+j)$.  Since $\tau(j)<p$, $\tau_j^p$ has a "new" descent at index $j$ in addition to all the descents of $\tau$, so $d_j(\tau_j^p)=d_j(\tau)+1$ and the two $(A,B)$ pairs are equal, as desired.  Similarly, if $\tau(j)>p$, then $mi(\tau,j+1,p)$ is computed using part (b) or (d) of Step 3, and so after this step we have $(A,B)=G(\tau,j)= (d_j(\tau),d_j(\tau)+j-1)$.  Since $\tau(j)>p$, $\tau_j^p$ has an ascent at index $j$ and thus $d_j(\tau)=d_j(\tau_j^p)$, so again the two $(A,B)$ pairs are equal.
\medskip \newline
The case of $q>p$ proceeds similarly, except that now $\chi(q>p)=1$, so we wish to show that each element of the first $(A,B)$ pair (computing the sequence $MIS(\tau_j^p,q)$) is greater by $1$ than the corresponding element in the second pair (computing the sequence $MIS(\tau,p)$).  The details are similar enough to the first case that we leave it to the reader to supply them. \hspace{\fill} \textbf{Q.E.D.}
\medskip \newline
To apply the lemma to our case, let $\tau=\sigma_{i+1}$, $p=\pi(i)$, and let $j=k_i$ (so that $\tau_j^p=\sigma_i$).  Finally, let $q=\pi(i-1)$.  Note that $m_i=MIS(\sigma_{i+1},\pi(i))(k_i)$, and that $m_{i-1}$ must be one of $MIS(\sigma_i,\pi(i-1))(1),...,MIS(\sigma_i,\pi(i-1))(k_i)$ because $\pi(i-1)$ must be inserted to the \emph{left} of $\pi(i)$.  We thus conclude, by the lemma, the following:  If $\pi(i-1)<\pi(i)$, then the only possible values of $m_{i-1}$ lie to the left of $m_i$ in $MIS(\sigma_{i+1},\pi(i))$, including $m_i$ itself; and if $\pi(i-1)>\pi(i)$, then the only possible values of $m_{i-1}$ are the values to the left of $m_i$ in $MIS(\sigma_{i+1},\pi(i))$ (again, including $m_i$ itself), each incremented by $1$.
\begin{description}
\item[Example 4.3.]  Let $\sigma_4=\theta=6152, \pi=437$, and $\sigma_1=\sigma=6143572$ (so that $\sigma_3=61572$ and $\sigma_2=613572$).  Then $k_3=4$ and $MIS(\sigma_4,7)=(3,2,4,\emph{1},0)$ (the italicized element being $m_3$).  By the lemma, we expect the first four elements of $MIS(\sigma_3,3)$--the sole candidates for $m_2$-- to be some permutation of $(3,2,4,1)$ (because $3<7$).  Indeed, the 4-initial segment of $MIS(\sigma_3,3)$ is $(2,3,\emph{1},4)$ (the italicized element denoting $m_2$, as $3$ is inserted at index $k_2=3$ to yield $\sigma_2$).  Again by the lemma, we expect the first three elements of $MIS(\sigma_2,4)$-- the candidates for $m_1$-- to be some permutation of $(3,4,2)$ (each of the first three elements of $MIS(\sigma_3,3)$ increased by $1$ because $4>3$) and indeed the $3$-initial segment of $MIS(\sigma_2,4)$ is $(3,4,2)$ itself.
\end{description}
We can now easily prove the following proposition:
\medskip \newline
\textbf{Proposition 4.1}: Let $S_i\subseteq[n]$ be the set of elements contained in the $k_i$-initial segment of $MIS(\sigma_{i+1},\pi(i))$ (for $i=1,...,a$) and let $T_i=S_i-d_i(\pi)=\{s-d_i(\pi)|s\in S_i\}$.  Then $T_1\subseteq...T_a\subseteq[b]$.
\medskip \newline
\textbf{Proof (of proposition)}:  By induction on the subscript of $T$, moving backwards from $a$ to $1$.  For $i=a$ this is simply an application of Proposition 3.1.  Suppose the proposition is true for $i=m+1,...,a$.  If $\pi(m)<\pi(m+1)$ then $d_m(\pi)=d_{m+1}(\pi)$ and $T_m\subseteq T_{m+1}$ iff $S_m\subseteq S_{m+1}$; if $\pi(m)>\pi(m+1)$ then $d_m(\pi)=d_{m+1}(\pi)+1$, and $T_m\subseteq T_{m+1}$ iff $S_m\subseteq \{s+1|s\in S_{m+1}\}$.  Both statements about the sets $S_m$ and $S_{m+1}$ are true by the lemma, as explained following the lemma and as illustrated in Example 4.3.  \hspace{\fill} \textbf{Q.E.D.}
\medskip \newline
Noting that $m_i\in S_i$ (as by definition, $m_i=MIS(\sigma_{i+1},\pi(i))(k_i)$), we immediately have:
\medskip \newline
\textbf{Corollary}: For all $i\in[a]$, $t_i=m_i-d_i(\pi)\in[b]$.
\medskip\newline
By this corollary, $0\leq t_i\leq b$ for all $i$, and thus $set((t_1,...,t_a))$ is a partition in $\mathcal{P}(b,a)$.  Hence $\Phi$ maps $\mathcal{S}(\theta,\pi)$ into $\mathcal{P}(b,a)$, as claimed.
\medskip \newline
It remains only to show that $\Phi$ is injective and surjective.  We do this by explaining how to find the unique $\sigma=\Phi^{-1}(\lambda)$ for any given partition $\lambda=(\lambda_1,...,\lambda_a)\in\mathcal{P}(b,a)$ ("unique", hence $\Phi$ is injective; "any", hence it is surjective).  The elements of $\lambda$ comprise one representative each from the sets $T_1,...,T_a$ defined in Proposition 4.1.   By that proposition, these sets form a nested chain.  For $i=a,...,1$, the choice of $m_i$ (and hence of $k_i$) determines both $t_i$ and the set $T_i$; specifically, the set $T_i$ contains precisely the first $k_i$ elements of $MIS(\sigma_{i+1},\pi(i))$ with $d_i(\pi)$ subtracted from each.  Thus the only way to ensure that $T_1\subseteq...\subseteq T_a$ is to choose $m_i$ to be the \emph{rightmost} element of $\{\lambda_i+d_i(\pi)|i=1,...,a\}$ which has not already been used in an earlier step, and thus $\Phi^{-1}(\lambda)$ is determined uniquely.  This completes the proof of Theorem 4.1. \hspace{\fill} \textbf{Q.E.D.}
\medskip \newline
We illustrate the method of determining $\Phi^{-1}(\lambda)$ using Example 4.1., now being performed in reverse.
\begin{description}
\item[Example 4.4.]  We compute the permutation $\sigma=\Phi^{-1}(\{0,3,4\})$ (where $\theta=5274$ and $\pi=631$ as in Example 4.1).  The values $0,3,$ and $4$ are the differences $m_i-d_i(\pi)$ for $i=1,2,3$.  For $i=3$, $d_3(\pi)$=0, and so $m_3$ must be $0,3,$ or $4$.  As $MIS(\theta,1)=(2,1,3,0,4)$, we must have $m_3=t_3=4$ (and hence $k_3=5$, $\sigma_3=52741$, and $T_3=\{0,1,2,3,4\}$) so that the remaining elements of $\lambda$ (0 and 3) are elements of $T_3$.
    \medskip \newline  For the next step (the $i=2$ step) we look at these remaining elements, each incremented by $d_2(\pi)=1$ to yield 1 and 4, and the $k_3$-initial segment of $MIS(\sigma_3,3)=(3,4,2,1,5,0)$.  The rightmost element in that segment among 1 and 4 is 1, at position 4, hence $m_2=1$ (and $t_2=0$), $k_2=4$, $\sigma_2=527341$, and $T_2=\{0,1,2,3\}$.
    \medskip \newline Finally, the remaining element of $\lambda$ is 3, now incremented by $d_1(\pi)=2$ to yield $m_1=5$.  The $k_2$-initial segment of $MIS(\sigma_2,6)=(4,3,2,5,6,1,0)$ indeed contains this value at the fourth position, and thus $k_1=4$ and $\sigma_1=\sigma=5276341$, as desired.
\end{description}
\section{Equidistribution over Inverse Descent Classes}\label{sec:equinv}
The \emph{inverse descent class} corresponding to a set $Q\subseteq[n]$ is the subset $S_Q\in S_n$ of all permutations of $[n]$ whose inverses (in the usual group-theoretic sense) have descent set $Q$.  There is a simple and well-known combinatorial description of inverse descent classes:  $k\in [n]$ is a descent of $\sigma^{-1}$ iff $k+1$ appears to the left of $k$ in $\sigma$.  Thus, if $Q=\{q_1,...,q_t\}$, then $S_Q$ is the set of shuffles of the complementary subsequences $\mathbf{q_0}=(1,...,q_1)$, $\mathbf{q_1}=(q_1+1,...,q_2)$,..., $\mathbf{q_t}=(q_t+1,...,n)$ such that none of these subsequences appears entirely to the right of any earlier subsequence.  We refer to a shuffle with this latter property as a \emph{well-mixed} shuffle.  It is generally easier to deal with all shuffles of $\mathbf{q_0},...,\mathbf{q_t}$ rather than only the well-mixed ones, and thus we focus not on the set $S_Q$ itself but rather on the larger set of permutations with inverses whose descent set is any subset of $Q$; by the combinatorial description above, this is precisely the set of all shuffles of $\mathbf{q_0},...,\mathbf{q_t}$.
\medskip\newline
Applying the theorem of Gessel and Garsia to this set yields an especially neat result because each subsequence $\mathbf{q_i} (i=0,...,t)$ is increasing, hence $maj(\mathbf{q_i})=0$.  Thus we have:
\begin{equation}\label{equ:idc_maj} \sum_{\{\sigma\in S_n|Des(\sigma^{-1}) \subseteq Q\}}q^{maj(\sigma)} = \left[
\begin{array}{c}
n \\ q_1,q_2-q_1...,q_t-q_{t-1},n-q_t
\end{array}
\right]
\end{equation}
It is shown in (\cite{Stanley}, Proposition 1.3.17) that the inversion number has the same generating function over the same set:
\begin{equation}\label{equ:idc_inv} \sum_{\{\sigma\in S_n|Des(\sigma^{-1}) \subseteq Q\}}q^{inv(\sigma)} = \left[
\begin{array}{c}
n \\ q_1,q_2-q_1...,q_t-q_{t-1},n-q_t
\end{array}
\right]
\end{equation}
The equidistribution of inversion number and major index over the set $\{\sigma\in S_n|Des(\sigma^{-1}) \subseteq Q\}$ follows immediately from these two equations, and their equidistribution over $S_Q$ itself follows from them as well by applying the inclusion-exclusion principle.
\medskip\newline
We conclude this paper by giving a direct bijective proof of these equidistribution results.  The proof addresses only the case of $|Q|=1$; specifically, we assume $Q=\{b\}$, with $n=a+b$, and (preserving the notation of the last section) $\theta=(1,...,b)$, $\pi=(b+1,...,n)$.  (For the general case of $Q=\{q_1,...,q_t\}$, the bijection is obtained by simply repeating the procedure described here $t$ times, where in the $i$-th round we assume $\theta$ to be any shuffle of $\mathbf{q_0},...,\mathbf{q_{i-1}}$ and $\pi=\mathbf{q_i}$).  Thus equations (\ref{equ:idc_maj}) and (\ref{equ:idc_inv}) become:
\begin{equation}\label{equ:idc_3} \sum_{\{\sigma\in S_n|Des(\sigma^{-1}) \subseteq Q\}}q^{maj(\sigma)} = \sum_{\{\sigma\in S_n|Des(\sigma^{-1}) \subseteq Q\}}q^{inv(\sigma)} =\left[
\begin{array}{c}
n \\ a
\end{array}
\right]
\end{equation}
Our approach will parallel the similar proof for all of $S_n$ given at the start of Section 3.  We will prove the generating function for inversion number by producing a bijection $\Psi: \mathcal{S}(\theta,\pi)\rightarrow \mathcal{P}(b,a)$ such that for $\tau\in\mathcal{S}(\theta,\pi)$, $inv(\tau)=|\Psi(\sigma)|$.  Then, utilizing the results of the last section, it will follow that $\Omega:=\Phi^{-1}\circ\Psi:S_Q\rightarrow S_Q$ is a bijection from $S_Q$ to itself which maps inversion number to major index, proving the equidistribution of these statistics over $S_Q$.
\medskip\newline
The bijection $\Psi$ is very simple.  Any shuffle $\tau$ of $\theta$ and $\pi$ is uniquely determined by the the weakly decreasing sequence $(t_1,...,t_a)$ where $t_i$ is the number of elements of $\theta$ to the right of $b+i$ in $\tau$.  Clearly, $0\leq t_i\leq b$ for all $i$, and conversely any sequence $(t_1,...,t_a)$ which is weakly decreasing with $0\leq t_i\leq b$ for all $i$ uniquely determines a shuffle $\tau$ of $\theta$ and $\pi$.  In this correspondence it is clear that $inv(\tau)=\sum_{i=1}^{a}t_i$.  Define $\Psi(\tau)=set(t_1,...,t_a)$.
\medskip\newline
To illustrate the bijection $\Omega:S_Q\rightarrow S_Q$ mapping inversion number to major index, consider the example of $n=7$, $b=4$, $\theta=(1,2,3,4)$, $\pi=(5,6,7)$, and $\tau=5126374$.  We have $\lambda:=\Psi(\tau)=\{4,2,1\}$.  We calculate $\sigma:=\Phi^{-1}(\lambda)$ using the method described at the end of the last section (which is especially simple here because $d_i(\pi)=0$ for all $i$ since $\pi$ is increasing).  As earlier, set $\sigma_4=\theta$, and the $i$-th insertion yields $\sigma_{4-i}$.  At each stage we italicize the element of the major increment sequence which is furthest to the right among the "unused" elements of $\lambda$; this determines both $m_i$ and $k_i$ at that stage.
\begin{center} \begin{tabular}{lll}
$MIS(\sigma_4,7)=(1,2,3,\emph{4},0)$ & $m_3=4,k_3=4$ & $\sigma_3=12374$\\
$MIS(\sigma_3,6)=(2,3,4,\emph{1},5,0)$ & $m_3=1,k_3=4$ & $\sigma_2=123674$\\
$MIS(\sigma_2,5)=(\emph{2},3,4,1,5,6,0)$ & $m_3=2,k_3=1$ & $\sigma_1=5123674$
\end{tabular} \end{center}
Thus we have $\Omega(5126374)=5123674$.
\medskip\newline
As a final remark, we note that $\Omega$ is not only a bijection on the set $\{\sigma\in S_n|Des(\sigma^{-1}) \subseteq Q\}$ but also on each individual inverse descent class contained in that set.  This is true because the shuffle of $\theta$ and $\pi$ which is not well-mixed (i.e., the shuffle $\sigma_0$ in which $\pi$ is appended to the right of $\theta$) is mapped to itself by $\Omega$:  $\Psi(\sigma_0)=(0,0,...,0)=\Phi^{-1}(\sigma_0)$.  When $\Omega$ is iterated multiple times for the case of $|Q|>1$, this fact remains true at each stage, and hence, for all $i\in[n-1]$, $i+1$ appears to the left of $i$ in $\Omega(\sigma)$ iff it does so in $\sigma$.  Thus the descents of $(\Omega(\sigma))^{-1}$ are the same as those of $\sigma^{-1}$, i.e., $\sigma$ and $\Omega(\sigma)$ are in the same inverse descent class.
\section{Appendix: Proving the L-G Algorithm}
We prove that the output permutation $\tau_1$ of the L-G Algorithm is $MIS(\sigma,r)$ by induction on the number of iterations.  We start by showing that the first iteration inserts the correct value of $mi(\sigma,n-1,r)$ to the left of $\tau_n$ and leaves the variables $A$ and $B$ assuming certain values which "set up" the next round; then we illustrate by induction that later iterations do the same.
\smallskip \newline
The first iteration implements part (a) or (d) of step 3 iff the last (rightmost) segment of $\sigma$ is a lesser segment, in which case (by step 2) $A=1, B=n-1$ at the start of this iteration, as in Algorithm L.  Thus part (a), implemented when the last segment is more than one letter long,  inserts $mi(\sigma,n-1,r)$ on the left of $\tau_n$ and yields $(A,B)=L(\sigma,n-2)$, as proven regarding Algorithm L.  Part (d) is implemented in the first iteration when the last segment contains only one letter, i.e., $\sigma(n-1)<r$ and $\sigma(n-2)>r$ (and hence $d_{n-2}=1$).  The insertion of $r$ between these two letters creates a new descent at index $n-1$, increasing the major index by $n-1=B$.  This is indeed the value that part (d) appends to the left of $\tau_n$.  $B$ is then decreased to $n-2$, and we have $d_{n-2}=1=A$, $(n-3)+d_{n-2}=n-2=B$, so this step yields $(A,B)=G(\sigma,n-2)$.
\smallskip \newline
The first iteration implements part (b) or (c) iff the last segment of $\sigma$ is greater, and step 2 sets up the initial values of $A=0$, $B=n-2$ as in Algorithm G.  Thus part (b), implemented when the last segment is more than one letter long,  inserts $mi(\sigma,n-1,r)$ on the left of $\tau_n$ and yields $(A,B)=G(\sigma,n-2)$, as proven regarding Algorithm G.  Part (c) is implemented in the first iteration when the last segment contains only one letter, i.e., $\sigma(n-1)>r$ and $\sigma(n-2)<r$ (and hence $d_{n-2}=0$).  The insertion of $r$ between these two letters creates no new descent, increasing the major index by $0=A$.  This is indeed the value that part (c) appends to the left of $\tau_n$.  When $A$ is increased to $1$ we have $1+d_{n-2}=1=A$, $(n-2)+d_{n-2}=n-2=B$, so this step yields $(A,B)=L(\sigma,n-2)$.
\smallskip \newline
Assume that the algorithm works in this manner for the first $j$ iterations, i.e., after $j$ iterations $\tau_{n-j}$ is identical to the last $j+1$ elements of $MIS(\sigma,r)$, and $(A,B)=L(\sigma,n-j-1)$ if the $j$-th iteration implemented part (a) or (c) in step 3, and $(A,B)=G(\sigma,n-j-1)$ if the $j$-th iteration implemented part (b) or (d) in step 3.  We show that the same remains true after the ($j+1$)-th iteration.
\smallskip \newline
Note that if the $j$-th iteration implemented part (a) or (c), then $\sigma(n-j-1)$ is in a lesser segment, so the $(j+1)$-th iteration must implement part (a) or (d).  By the inductive assumption, in this case we have $(A,B)=L(\sigma,n-j-1)$ after the $j$-th iteration, so if the $(j+1)$-th iteration implements part (a) then it appends the correct value to the left of $\tau_{n-j}$ and (unless this iteration is the last one, in which case the algorithm ends) yields $(A,B)=L(\sigma,n-j-2)$ by the proof of Proposition 3.2. If the $(j+1)$-th iteration implements part (d), i.e., $\sigma(n-j-1)<r$ and $\sigma(n-j-2)>r$ (so that $d_{n-j-2}=d_{n-j-1}+1$), then inserting $r$ between these two letters creates a new descent at index $n-j-1$, increasing the major index by $(n-j-1)+d_{n-j-1}$.  By assumption, $(A,B)=L(\sigma,n-j-1)=(1+d_{n-j-1},(n-j-1)+d_{n-j-1})$, so the increase in major index is exactly $B$, hence part (d) works correctly.  $B$ is then decreased by 1 to yield $(A,B)= (1+d_{n-j-1},(n-j-1)+d_{n-j-1}-1)= (d_{n-j-2},(n-j-3)+d_{n-j-2})=G(\sigma,n-j-2)$, completing the induction.
\smallskip \newline
Similarly, if the $j$-th iteration implemented part (b) or (d), then $\sigma(n-j-1)$ is in a greater segment, so the $(j+1)$-th iteration must implement part (b) or (c).  By the inductive assumption, in this case we have $(A,B)=G(\sigma,n-j-1)$ after the $j$-th iteration, so if the $(j+1)$-th iteration implements part (b) then it appends the correct value to the left of $\tau_{n-j}$ and (unless this iteration is the last one, in which case the algorithm ends) yields $(A,B)=G(\sigma,n-j-2)$ by the proof of Algorithm G above. If the $(j+1)$-th iteration implements part (c), i.e., $\sigma(n-j-1)>r$ and $\sigma(n-j-2)<r$ (so that $d_{n-j-2}=d_{n-j-1}$), then inserting $r$ between these two letters creates no new descent at all, so the major index increases by $d_{n-j-1}$.  By assumption, $(A,B)=G(\sigma,n-j-1)=(d_{n-j-1},(n-j-2)+d_{n-j-1})$, so the increase in major index is exactly $A$, hence part (c) works correctly.  $A$ is then increased by 1 to yield $(A,B)= (1+d_{n-j-1},(n-j-2)+d_{n-j-1})= (1+d_{n-j-2},(n-j-2)+d_{n-j-2})=L(\sigma,n-j-2)$, completing the induction.
\begin{flushright} \textbf{Q.E.D.} \end{flushright}

\end{document}